\documentclass[twoside,11pt]{amsart}
\usepackage{amsmath,latexsym,amssymb,mathptm, enumerate,verbatim}
\usepackage{diagrams}
\def\proof{\noindent{\bf{Proof.} }}
\def\sqr#1#2{{\vcenter{\hrule height.#2pt
        \hbox{\vrule width.#2pt height#1pt \kern#1pt
                \vrule width.#2pt}
        \hrule height.#2pt}}}

\def\tratto{\mbox{\rule{2mm}{.2mm}$\;\!$}}

\newtheorem{theorem}{Theorem}[section]

\newtheorem{corollary}[theorem]{Corollary}
\newtheorem{lemma}[theorem]{Lemma}
\newtheorem{proposition}[theorem]{Proposition}

 \theoremstyle{definition}

\newtheorem{definition}[theorem]{Definition}
\newtheorem{discussion}[theorem]{Discussion}

\newtheorem{remark}[theorem]{Remark}

\newtheorem{example}[theorem]{Example}
\newtheorem{question}[theorem]{Question}

\newcommand{\s}{\bigskip}

\newcommand{\del}{\ensuremath{\partial}}
\newcommand{\D}{\ensuremath{\Delta}}
\newcommand{\sd}{\ensuremath{\delta}}

\newcommand{\m}{\mathfrak{m}}
\newcommand{\la}{\lambda}
\newcommand{\lr}{\longrightarrow}
\newcommand{\ol}[1]{\ensuremath{\overline{#1}}}
\newcommand{\ul}[1]{\ensuremath{\underline{#1}}}

\numberwithin{equation}{section}

\newcommand{\ds}{\displaystyle}

\begin{document}
\baselineskip=16pt

\title{What is a system of parameters?}
\author[L. Fouli]{Louiza Fouli}
\address{Department of Mathematical Sciences, New Mexico State University, Las Cruces, NM 88003}
\email{lfouli@math.nmsu.edu}
\urladdr{http://www.math.nmsu.edu/\textasciitilde lfouli}

\author[C. Huneke]{Craig Huneke}
\address{Department of Mathematics, University of Kansas, Lawrence, KS 66045}
\email{huneke@math.ku.edu}
\urladdr{http://www.math.ku.edu/\textasciitilde huneke}

\thanks{2010 {\it Mathematics Subject Classification.} 13A35, 13C40, 13D45}

\thanks{The first author was partially supported by the NSF-AWM Mentoring Travel Grant, grant DMS-0839954 and second was partially supported by the National Science Foundation, grant
DMS-0756853}

\keywords{system of parameters, local cohomology, limit closure, tight closure, plus closure}

\vspace{-0.1in}

\begin{abstract} In this paper we discuss various refinements and generalizations of a theorem of Sankar Dutta and
Paul Roberts. Their theorem gives a criterion for $d$ elements in a $d$-dimensional Noetherian Cohen-Macaulay local ring to
be a system of parameters, i.e., to have height $d$. We chiefly remove the assumption that the ring be
Cohen-Macaulay and discuss similar theorems.
\end{abstract}

\vspace{-0.2in}

\maketitle

\section{Introduction}
One of the fundamental concepts in commutative algebra is
that of \it height\rm. Let $(R,\m)$ be a Noetherian local ring of
dimension $d$. Every ideal $I$ of $R$ is generated by
$d$ elements up to radical. The question of determining the height of ideals in general is then
the same as that of determining the height of an ideal generated by
$d$ elements. An interesting paper of Sankar Dutta and Paul Roberts \cite{DR}
gives a criterion for $d$ elements in a Cohen-Macaulay local ring to be a system of parameters, i.e., to
have maximal height. We recall their theorem: let $(R,\m)$ be a $d$-dimensional
Cohen-Macaulay local ring, and let $x_1,\ldots, x_d$ be a system of parameters.
Suppose that $(y_1,\ldots, y_d)\subset (x_1,\ldots, x_d)$. Write
$y_i = \sum \limits_j a_{ij}x_j$ (since this notation occurs throughout this paper, we abbreviate
it by writing $ (\ul{y}) \overset{A}{\subset} (\ul{x})$), and let $\det A$ be the determinant of the matrix
whose coefficients are the $a_{ij}$. Since $\det A\cdot (x_1,\ldots, x_d)\subset
(y_1,\ldots, y_d)$, multiplication by $\det A$ induces a well-defined map
from $R/(x_1,\ldots, x_d)$ to $R/(y_1,\ldots, y_d)$. The theorem of
Dutta and Roberts states that this map is injective if and only if $y_1,\ldots,y_d$
form a system of parameters, i.e., if and only if the ideal $(y_1,\ldots, y_d)$ has
height $d$. The fact that the map is injective if both $y_1,\ldots,y_d$ and $x_1,\ldots, x_d$ are regular sequences
was well-known; it is the converse that was new.

In the 1970s, Hochster introduced the problem of understanding
what constraints parameters in a Noetherian local ring must satisfy. This means understanding what
equations \it cannot \rm be satisfied by a system of parameters in a Noetherian local ring. The most famous
such constraint is the monomial conjecture of Hochster.  The theorem of Dutta and Roberts can
be thought of as giving the opposite: it gives a relation which can \it only \rm be satisfied by
parameters. 
However, the Dutta-Roberts result leads to many interesting questions, some of which turn out to be closely related,
in various guises, to the monomial conjecture. We will mention some of these questions in the last section. A first and obvious question is whether
or not the same result holds without the Cohen-Macaulay assumption. We have been able to
answer this for 1-dimensional rings (Theorem~\ref{1-dim DR general} and Example~\ref{ex Heitmann}). Even this case is not trivial.
We have no counterexample in higher dimension.

Perhaps an even more basic question is
whether or not the assumption of the theorem is independent of the matrix of coefficients
chosen. For example, a famous special case of the theorem is the case in which $R$ is a polynomial ring, $k[x_1,..,x_d]$ over a
field of characteristic $0$, and $y_1,\ldots,y_d$ are homogeneous of degree $n_1,\ldots,n_d$ respectively. In this
case Euler's formula gives that $y_i =\ds{ \frac{1}{n_i}\sum \limits_{1\leq j\leq d}\frac{\partial{y_i}}{\partial{x_j}}x_j}$, and letting $\Delta$ be the
determinant of this particular matrix of coefficients, one has that $\Delta\notin (y_1,\ldots,y_d)$ if and only
if the height of $(y_1,\ldots,y_d)$ is $d$ \cite{SS}. In the language above, the map from $R/(x_1,\ldots,x_d)$ to
$R/(y_1,\ldots,y_d)$ is either injective or zero. Is this true for an arbitrary matrix of coefficients? Strooker \cite{St}
proved that independence is guaranteed if one changes the statement of the
theorem to apply to the map induced by multiplication by the determinant from
$R/(x_1,\ldots, x_d)^{\lim}$ to $R/(y_1,\ldots, y_d)^{\lim}$ (part of the statement is that
this is a well-defined homomorphism). Here $R/(x_1,\ldots, x_d)^{\lim}$ (respectively
 $R/(y_1,\ldots, y_d)^{\lim}$) denotes the image of $R/(x_1,\ldots, x_d)$ (respectively
 $R/(y_1,\ldots, y_d)$) under the identification of $H^d_{\m}(R)$ with the direct limit
of $R/(x_1^n,\ldots, x_d^n)$ (respectively the image of  $R/(y_1,\ldots, y_d)$ in
$H^d_{(y_1,\ldots,y_d)}(R)$ under the identification of this module with
the direct limit of $R/(y_1^n,\ldots, y_d^n)$). Thus the question of independence of
the choice of matrix naturally leads to a consideration of other maps, and in particular brings into the picture
local cohomology.

We begin this paper by proving various maps are always injective if both  the sequence $\ul{x}:=x_1, \ldots, x_d$ and the sequence
$\ul{y}=y_1, \ldots, y_d$ form systems of parameters. Included in this list are the map above from
$R/(\ul{x})^{\lim}$ to $R/(\ul{y})^{\lim}$ induced by multiplication by
$\det A$ (this is now independent of the matrix of coefficients $(a_{ij})$ chosen),
a natural map from $H^d_{(\ul{x})}(R)$ to $H^d_{(\ul{y})}(R)$, and in positive prime characteristic,
the map from $R/(\ul{x})^*$ to $R/(\ul{y})^*$  induced by multiplication by $\det A$,
where the $^*$ refers to tight closure. One can also replace the tight closure by plus closure
(they are the same for parameters by \cite{Sm}, but are not necessarily the same for
other ideals).

We next ask whether the obvious generalization of the theorem of Dutta and Roberts holds in
the cases listed in the paragraph above: if the maps are injective, must the $\ul{y}$ form
a system of parameters? We are able to prove that there is a fixed power of $\m$, independent of
 $\ul{x}$ and $\ul{y}$, such that if all the $x_i$ are in  this fixed power, then the natural
generalization of the theorem of Dutta and Roberts holds in the first two maps in the paragraph above (Corollary~\ref{lim inj iff par in m^l}).
Moreover, if the ring is analytically irreducible, then we obtain a full generalization (Proposition~\ref{Hart-Lich iff}). This is
mainly due to the vanishing theorem of Hartshorne and Lichtenbaum, which can be thought of as
giving a necessary and sufficient condition for $d$ elements $\ul{y}$ in a complete local domain of dimension
$d$ to be a system of parameters, namely that $H^d_{(\ul{y})}(R) \ne 0$. In the cases of tight closure
and plus closure, we do not know whether injectivity on the maps forces the $\ul{y}$ to be a
system of parameters. We have been unable to prove this or give a counterexample.

As we mentioned above, a full generalization of the Dutta-Roberts theorem escapes us at the moment. Note, however,
that the well-known direction of their theorem is false in the non Cohen-Macaulay case, unlike for the other
maps discussed above. In fact we prove that for all systems of parameters $\ul{x}$ and all
$ (\ul{y}) \overset{A}{\subset} (\ul{x})$ the induced map $R/(\ul{x}) \overset{\cdot \det A}\to R/(\ul{y})$  is injective  if and only if the ring is Cohen-Macaulay.

We do give the best possible answer for the limit closures. Our final result proves the following:
Let $(R,\m)$ be a local ring of equicharacteristic or of dimension at most three. Set $d = \dim R$.  There exists an integer $\ell$ with the following property:
whenever $\ul{x}=x_{1}, \ldots, x_d$ is a system of parameters with $(\ul{x}) \subset \m^{\ell}$ and  $\ul{y}=y_1, \ldots, y_d$
a sequence of elements such that $(\ul{y}) \stackrel{A}{\subset} (\ul{x})$, then $\ul{y}$ is a system of parameters if and only if the map $R/{(\ul{x})^{\rm lim}} \stackrel{\cdot \det A}{\lr} R/{(\ul{y})^{\rm lim}}$ is injective.
We need $R$ to be equicharacteristic or of dimension at most three in order to apply the monomial conjecture, which  is known to hold in these cases
\cite{He}. We give an example to show that this result is not true without the parameters being ``deep" enough inside the maximal ideal (Example~\ref{highpower}).

\section{Preliminaries and Basic Results}

Let $\ul{x}=x_1, \ldots, x_n$ be a sequence of elements in a Noetherian ring $R$. Let $(\ul{x})$ denote the ideal generated
 by $x_1, \ldots, x_n$. We recall here that the limit closure of $(\ul{x})$ is given by
 $$(\ul{x})^{\rm lim}=\bigcup_{t \geq 1} (x_1^t,\ldots, x_n^{t}):x_1^{t-1} \cdots x_n^{t-1}.$$

Before stating some of the conditions we will be studying, we need to recall some well-known definitions.

\begin{definition} Let $R$ be a Noetherian local ring of positive characteristic $p$. An element $x$ is in the tight closure
of an ideal $I$ if there exists an element $c$, not in any minimal prime of $R$, such that for all large $q = p^e$,
$cx^q\in I^{[q]}$, where $I^{[q]}$ is the ideal generated by all $f^q$ for $f\in I$. \end{definition}

\begin{definition}  Let $R$ be a Noetherian domain. The absolute integral closure of $R$, denoted $R^+$, is the
integral closure of $R$ in a fixed algebraic closure of the fraction field of $R$. \end{definition}

\begin{remark}\label{equivalences}
{\rm Let $(R,\m)$ be a Noetherian local ring of dimension $d$. Let $\ul{x}=x_1,\ldots, x_d$ be a system of parameters and let
 $\ul{y}=y_1, \ldots, y_d$ be a sequence of elements in $R$ such that $(\ul{y}) \subset (\ul{x})$.
 Let $A=(a_{ij})$ be a matrix such that $y_i=\sum \limits_{j=1}^{d} a_{ij}x_j$ and let $\det A$ denote the determinant of $A$. We will use the notation $(\ul{y}) \stackrel{A}{\subset} (\ul{x})$ to denote that $A$ is a matrix of coefficients relating $\ul{y}$ to $\ul{x}$. Consider the following
 statements, where in parts (3) and (4) we assume that $R$ has positive characteristic $p$, and in case (4) that $R$ is a domain:

\begin{enumerate}[{\bfseries (1)}]

\item The map $R/{(\ul{x})^{\rm lim}} \overset{{\cdot \det A}}{\longrightarrow} R/{(\ul{y})^{\rm lim}}$ is injective.

\item The map $H_{(\ul{x})}^{d}(R) \longrightarrow H_{(\ul{y})}^{d}(R)$ is injective.

\item The map $R/{(\ul{x})^{*}} \overset{{\cdot \det A}}{\longrightarrow} R/{(\ul{y})^{*}}$ is injective.

\item The map $H_{(\ul{x})}^{d}(R^{+}) \longrightarrow H_{(\ul{y})}^{d}(R^{+})$ is injective.

\item The map $R/(\ul{x}) \stackrel{\cdot \det {A}}{\lr} R/(\ul{y})$ is injective.

\item $\ul{y}$ is a system of parameters.
\end{enumerate}}
\end{remark}

The maps on local cohomology are not yet defined, but will be defined in a way compatible with the maps on the lim closures.
We wish to determine the relationship between these statements. In this section we will prove that the maps in parts
(1) and (2) are injective if $\ul{y}$ is a system of parameters. In a later section we will prove that the maps
in (3) and (4) are injective if $\ul{y}$ is a system of parameters. In contrast, the map in (5) does not necessarily have to be
injective when (6) holds.

The fact that the map in (1) is well-defined (and does not depend on the specific matrix $A$) was proved by Strooker
in \cite{St}. We will prove this as well as a by-product of the following lemma, which we need to define the map
in (2).

\begin{lemma}\label{homotopyfrombelow}
Let $R$ be a commutative ring with identity. Let $(G_{\bullet},\Delta)$ and $(F_{\bullet},\partial)$ be complexes of $R$-modules. 
Let $\alpha: G_{\bullet} \rightarrow F_{\bullet}$ be a map of complexes. Suppose that for all $0\leq j \leq i-1$ there
 exist maps $\delta_{j}:G_j \rightarrow F_{j+1}$ such that $\alpha_j=\del_{j+1}\sd_j+\sd_{j-1}\D_j$. We further
 assume that there exists a complex of free $R$-modules $({G'}_{\bullet},\D')$ and a map of complexes
 $\beta: {G'}_{\bullet} \rightarrow G_{\bullet}$ such that $I_{1}(\beta_{i}) H_{i}(F_{\bullet})=0$. Then
 we may extend the homotopy to $\sd_i:G'_{i} \rightarrow F_{i+1}$ for the map of complexes $\alpha \circ \beta : {G'}_{\bullet} \rightarrow F_{\bullet}$.
\end{lemma}

\proof
We have the following commutative diagram:

\begin{diagram}[width=2em]
 &\cdots &\rTo &F_{i+1}&\rTo ^{\del_{i+1}} &F_i &\rTo^{\del_{i}} &F_{i-1} &\rTo & \cdots&\rTo &F_1 &\rTo^{\del_1} &F_0 &\rTo &0\\
 &&&\uTo^{\alpha_{i+1}} &&\uTo^{\alpha_{i}} &\luTo(2,4)\;\;\;\;\;\;\;\;\;\;{\sd_{i-1}'}\luTo(2,2)>\; \; \; \; \;\; \;{\sd_{i-1}} &\uTo_{\alpha_{i-1}} &&&&\uTo^{\alpha_1}&\luTo(2,4)\;\;\;\;\sd_{0}' \;\;\;\;\luTo(2,2)\;\;\;\;\sd_{0}&\uTo_{\alpha_0}&&\\
 &\cdots &\rTo &G_{i+1}&\rTo^{\D_{i+1}} &G_i &\rTo^{ \; \; \; \; \; \; \; \D_{i}} &G_{i-1} &\rTo & \cdots&\rTo &G_1 &\rTo^{\;\;\;\;\; \D_1} &G_0 &\rTo &0\\
 &&&\uTo^{\beta_{i+1}} &&\uTo^{\beta_{i}} & &\hspace{0.3cm}\uTo_{\beta_{i-1}} &&&&\uTo^{\beta_1}&&\uTo_{\beta_0}&&\\
 &\cdots &\rTo &G_{i+1}'&\rTo ^{\D_{i+1}'} &G_i' &\rTo^{\D_{i}'} &G_{i-1}' &\rTo & \cdots&\rTo &G_1' &\rTo^{\D_1'} &G_0' &\rTo &0\\
 \end{diagram}

First we begin by showing the there is a homotopy for the map of complexes $\alpha \circ \beta : {G'}_{\bullet} \rightarrow F_{\bullet}$ up
to level $i-1$. For $0 \leq j \leq i-1$ define $\sd'_j=\sd_j \circ \beta_j$.
Then $\sd_j':{G'}_j \rightarrow F_{j+1}$ and $\alpha_j \circ \beta_j=(\del_{j+1}\sd_j+\sd_{j-1}\D_j)
\beta_j=\del_{j+1}\sd_j\beta_j+\sd_{j-1}\D_j\beta_j=\del_{j+1}\sd_j'+\sd_{j-1}\beta_{j-1}\D_j'=\del_{j+1}\sd_j'+\sd_{j-1}'\D_j'$.

We now claim that ${\rm Im} (\alpha_i-\sd_{i-1}\D_i) \subset \ker \del_i$. Indeed, $$\del_i(\alpha_i-\sd_{i-1}\D_i)=\del_i \alpha_i-\del_i\sd_{i-1}\D_i=\alpha_{i-1}\D_i -\del_i\sd_{i-1}\D_i=\del_{i}\sd_{i-1}\D_i+\sd_{i-2}\D_{i-1}\D_i-\del_i\sd_{i-1}\D_i=0.$$

Finally we show that ${\rm Im} (\alpha_i\beta_i-\sd_{i-1}' \D_i') \subset {\rm Im}(\del_{i+1})$. To see this notice that ${\rm Im} (\alpha_i\beta_i-\sd_{i-1}' \D_i')={\rm Im} (\alpha_i\beta_i-\sd_{i-1}\beta_{i-1} \D_i')={\rm Im} (\alpha_i\beta_i-\sd_{i-1}\D_i\beta_{i}) \subset {\rm Im} (\alpha_i-\sd_{i-1}\D_i)I_{1}(\beta_{i}) \subset {\rm Im}(\del_{i+1})$, since we have $I_{1}(\beta_{i}) H_{i}(F_{\bullet})=0$.
Therefore we may lift to a map $\sd_{i}':G_i' \rightarrow F_{i+1}$ such that $\alpha_i\beta_i-\sd_{i-1}' \D_i'=\sd_{i}'\del_{i+1}$. \qed

\s
As an immediate corollary, we get the following effective result which explains how the different choices of matrics relate to the 
multiplication maps of the determinants of the matrices.

\begin{corollary}\label{det-cor}
Let $R$ be a commutative ring with identity, $d=\dim R$ and let $(y_1, \ldots, y_d) \subset (x_1, \ldots, x_d)$. Suppose that there
exist two matrices $A=(a_{ij})$ and $B=(b_{ij})$ such that $y_i=\sum \limits_{j=1}^{d} a_{ij}x_j=\sum \limits_{j=1}^{d} b_{ij}x_j$. 
Then $(y_1 \cdots y_d)^d(\det A- \det B) \in (y_1^{d+1}, \ldots, y_d^{d+1})$.
\end{corollary}

\proof Let $\ul{y}=y_1, \ldots, y_d$ and $\ul{x}=x_1, \ldots, x_d$.
Consider the map of the Koszul complexes $K_{\bullet}(\ul{y};R) \rightarrow K_{\bullet}(\ul{x};R)$:
\begin{diagram}[width=2em]
\cdots &R^{d}&\rTo^{(x_1, \cdots, x_d)}&R&\rTo&R/(\ul{x})&\rTo &0\\
&\uTo^{A}&& \parallel&&\uTo\\
\cdots &R^{d}&\rTo^{(y_1, \cdots, y_d)}&R&\rTo&R/(\ul{y})&\rTo &0\\
\end{diagram}

The map induced by $A$ extends to the whole Koszul complex by using exteria powers of $A$, $\wedge^{i}A:\wedge^{i}(R^{d}) \rightarrow \wedge^{i}(R^{d})$.
Similarly we obtain the map $K_{\bullet}(\ul{y};R) \stackrel{\wedge^{\bullet}}{\rightarrow} K_{\bullet}(\ul{x};R)$.
Note that we have the following diagram:
\begin{diagram}[width=2em]
\cdots &R^{d}&\rTo^{(x_1, \cdots, x_d)}&R&\rTo&R/(\ul{x})&\rTo &0\\
&\uTo^{A-B}&& \uTo^{0}&&\uTo\\
\cdots &R^{d}&\rTo^{(y_1, \cdots, y_d)}&R&\rTo&R/(\ul{y})&\rTo &0\\
\end{diagram}

On the other hand $H_{\bullet}(K_{\bullet}(\ul{x};R))$ is annihilated by $(\ul{x})$ and therefore it is also annihilated by $(\ul{y})$. We will begin defining a homotopy for the map of complexes $\wedge^{i}(A)-\wedge^{i}(B): K_{\bullet}(\ul{y};R)\rightarrow K_{\bullet}(\ul{x};R)$ by taking $\sd_0=0$. We compose with $K_{\bullet}(y_1^2,\ldots, y_d^2;R) \rightarrow K_{\bullet}(\ul{y};R)$ induced by $\wedge^{\bullet} \small{\left(\begin{array}{ccc}y_1 &&0\\&\ddots& \\ 0&&y_d\end{array}\right)}$. Therefore by Lemma~\ref{homotopyfrombelow} we may extend the homotopy to $\sd_1$ for $K_{\bullet}(y_1^2,\ldots, y_d^2;R) \rightarrow K_{\bullet}(\ul{x};R)$. Repeating the process we obtain:

\small{\begin{diagram}[width=0.5em]
K_{\bullet}(\ul{x};R):&0&\rTo&R&\rTo&R^d&\rTo&\cdots &\rTo&R^{d}&\rTo&R&\rTo&R/(\ul{x})&\rTo &0\\
&&&\uTo^{\;\;\;(y_1\cdots y_d)^{d-1}(\det A-\det B)}&\luTo(2,2)\sd_{d-1}&\uTo&&&\luTo(2,2)\; \; \; \; \; \; \; {\sd_1}&\uTo&\luTo(2,2)\; \; \; \; \; \; {\sd_{0}} & \uTo_{0}&&\uTo&&\\
K_{\bullet}(\ul{y}^{[d]};R):&0&\rTo&R&\rTo^{(\ul{y}^{[d]})\;\;\;\;\;\;\;} &R^d&\rTo&\cdots&\rTo&R^{d}&\rTo&R&\rTo&R/(\ul{y}^{[d]})&\rTo &0\\
\end{diagram}}

Hence $(y_1 \cdots y_d)^{d-1}(\det A-\det B) \in (y_1^d,\ldots, y_d^{d})$. Finally multiplying by $(y_1 \cdots y_d)$
yields $$(y_1 \cdots y_d)^{d}(\det A-\det B) \in y_1 \cdots y_d(y_1^d,\ldots, y_d^{d}) \subset (y_1^{d+1},\ldots, y_d^{d+1}). \;\;\;\;\qed$$

\s

We can now give a quick proof for the map in  part~(1) from Remark~\ref{equivalences} to be well defined. This is due to Strooker. See \cite[5.1.14--5.1.17]{St}.

\begin{remark}\label{DeltaMap}
{\rm Let $R$ be a commutative ring with identity, $d=\dim R$ and let $\ul{x}=x_1, \ldots, x_d$ and $\ul{y}=y_1, \ldots, y_d$ be sequences of elements such that $(\ul{y}) \stackrel{A}{\subset}  (\ul{x})$. 
Then the map $R/{(\ul{x})^{\rm lim}} \overset{{\cdot \det A}}{\longrightarrow} R/{(\ul{y})^{\rm lim}}$ induced by multiplication by $\det A$ is well defined .
Moreover, this homomorphism does not depend on the choice of $A$; if $(\ul{y}) \stackrel{B}{\subset} (\ul{x})$, then multiplication by
$\det B$ induces the same map.}
\end{remark}

\proof
Let $r \in (\ul{x})^{\rm lim}$. Then there exists a positive integer $t$ such that $r \cdot x_1^{t-1} \cdots x_n^{t-1} \in (x_1^t,\ldots, x_n^{t}):=(\ul{x})^{[t]}$. Since $(\ul{y}) \subset (\ul{x})$ then there exists a positive integer $s$ such that $(\ul{y})^{[s]} \subset (\ul{x})^{[t]}$. Let $D$ be a matrix such that $(\ul{y})^{[s]}=D (\ul{x})^{[t]}$.
Let $C_1$ be the diagonal matrix with entries $y_i^{s-1}$ and $C_2$ be the diagonal matrix with entries $x_i^{t-1}$. Notice that we have the following inclusions:
\begin{align*}
&(\ul{y})^{[s]} \overset{C_1}{\subset} (\ul{y}) \overset{A}{\subset} (\ul{x}) \\
&(\ul{y})^{[s]} \overset{D}{\subset} (\ul{x})^{[t]} \overset{C_2}{\subset} (\ul{x})
\end{align*}

Let $C=AC_1$ and $E=C_2D$. Then by Corollary~\ref{det-cor}
$(y_1 \cdots y_d)^{sd}(\det C- \det E) \in (\ul{y})^{[sd+s]}$. Hence $$ \det A (y_1 \cdots y_d)^{s-1}(y_1 \cdots y_d)^{sd}-\det D(y_1 \cdots y_d)^{sd}(x_1 \cdots x_d)^{t-1} \in (\ul{y})^{[sd+s]}$$ and thus $r\det A (y_1 \cdots y_d)^{s-1}(y_1 \cdots y_d)^{sd}-r
\det D (y_1 \cdots y_d)^{sd}(x_1 \cdots x_d)^{t-1} \in (\ul{y})^{[sd+s]}$. Notice, however that $r\det D (y_1 \cdots y_d)^{sd}(x_1 \cdots x_d)^{t-1} \in \det D(y_1 \cdots y_d)^{sd} (\ul{x})^{[t]}$, since
$r \in (\ul{x})^{\rm lim}$. Also, $$\det D (y_1 \cdots y_d)^{sd}(\ul{x})^{[t]} \subset  (y_1 \cdots y_d)^{sd}(y_1 \cdots y_d)^{[s]} \subset (\ul{y})^{[sd+s]}.$$ 
Thus $r\det A(y_1 \cdots y_d)^{s-1}(y_1 \cdots y_d)^{sd}=r \det A(y_1 \cdots y_d)^{sd+s-1} \in (\ul{y})^{[sd+s]}$. Therefore, $r \in (\ul{y})^{\rm lim}$ and the map is then well defined.

The fact that the map does not depend on the choice of $A$ follows directly from Corollary~\ref{det-cor}; this says that the difference between
the two determinants lies in $(\ul{y})^{\rm lim}$.
\qed

\s

Next we define the maps from part~(2) and equivalently part~(4) from Remark~\ref{equivalences}.

\begin{discussion}\label{Phi Disc} {\rm Let $R$ be a Noetherian local ring, $d=\dim R$ and let $\ul{x}=x_1, \ldots, x_d$ and $\ul{y}=y_1, \ldots, y_d$ be sequences of elements such that $(\ul{y}) \stackrel{A}{\subset}  (\ul{x})$.  We will define a homomorphism $\Phi: H_{(\ul{x})}^{d}(R) \longrightarrow H_{(\ul{y})}^{d}(R)$ such that
the following diagram commutes:

\begin{diagram}[width=2em]
0 &\rTo& R/{(\ul{x})^{\rm lim}} &\rTo^{\alpha} & H_{(\ul{x})}^{d}(R) \\
&& \dTo^{\beta=\cdot \det A} & & \dTo^{\Phi} \\
0 &\rTo& R/{(\ul{y})^{\rm lim}} &\rTo^{\gamma} & H_{(\ul{y})}^{d}(R) \\
\end{diagram}

In other words, when  
we restrict $\Phi$ to the natural image of $R/(\ul{x})$ in  $H_{(\ul{x})}^{d}(R)$, the map goes to $\det A$ times the image of  $R/(\ul{y})$ in  $H_{(\ul{y})}^{d}(R)$,
i.e., is compatible with the map given in part~(1) of Remark~\ref{equivalences}. To define $\Phi$, let
$u \in H_{(\ul{x})}^{d}(R)$. Then $u=\ds{[\frac{r}{x_1^{n} \cdots x_d^{n}}]}$ for some $n$ and some $r \in R$.
Since $(\ul{y}) \subset (\ul{x})$ there exists a matrix $B$ and a positive integer $s$ such that
$(\ul{y})^{[s]}=B(\ul{x})^{[n]}$. We define $\Phi: H_{(\ul{x})}^{d}(R) \longrightarrow H_{(\ul{y})}^{d}(R)$ 
by $\Phi(u)=\ds{[\frac{r \det B}{y_1^{s} \cdots y_d^{s}}]}$.}
\end{discussion}

\begin{proposition}\label{Hmap}
Let $R$ be a Noetherian local ring, $d=\dim R$ and let $\ul{x}=x_1, \ldots, x_d$ and $\ul{y}=y_1, \ldots, y_d$ be sequences of elements such that $(\ul{y}) \stackrel{A}{\subset}  (\ul{x})$. Then $\Phi$ is a homomorphism from $H_{(\ul{x})}^{d}(R) \longrightarrow H_{(\ul{y})}^{d}(R)$,
and $\Phi(\ds{[\frac{1}{x_1\cdots x_d}]}) = \ds{[\frac {\det A}{y_1\cdots y_d}]}$.
\end{proposition}

\proof
We adopt the notation of the above discussion.
Suppose that $\ds{[\frac{r_1}{x_1^{n} \cdots x_d^{n}}]}=\ds{[\frac{r_2}{x_1^{n} \cdots x_d^{n}}]}$, where $r_1, r_2 \in R$. 
Then there exists a positive integer $t$ such that $(r_1-r_2)(x_1 \cdots x_d)^t \in (\ul{x})^{[t+n]}$. We may assume that $t >n$.
There exists a matrix $C$ and a positive integer $l$ such that $(\ul{y})^{[l]}=C (\ul{x})^{[t+n]}$. We have the following diagram:

\begin{align*}
&(\ul{y})^{[l]} \overset{C}{\subset} (\ul{x})^{[t+n]} \overset{D_{(\ul{x})}}{\subset} (\ul{x})^{[n]}, \\
&(\ul{y})^{[l]} \overset{D_{(\ul{y})}}{\subset} (\ul{y})^{[s]} \overset{B}{\subset} (\ul{x})^{[n]},
\end{align*}
where $D_{(\ul{x})}$ is the diagonal matrix with entries $x_i^{t}$ and $D_{(\ul{y})}$ is the diagonal matrix with entries $y_i^{l-s}$.

Therefore by Corollary~\ref{det-cor} $(y_1 \cdots y_d)^{dl}(\det B (y_1 \cdots y_d)^{l-s}-(x_1 \cdots x_d)^t \det C) \in (\ul{y})^{[dl+l]}$. 
Multiplying by $r_1-r_2$ we obtain $(r_1-r_2) \det B (y_1 \cdots y_d)^{dl+l-s}-(r_1-r_2)(y_1 \cdots y_d)^{dl}(x_1 \cdots x_d)^t
\det C \in (\ul{y})^{[dl+l]}$. Notice that since $(r_1-r_2)(x_1 \cdots x_d)^t \in (\ul{x})^{[t+n]}$
then $(r_1-r_2)(y_1 \cdots y_d)^{dl}(x_1 \cdots x_d)^t \det C \in (y_1 \cdots y_d)^{dl} (\ul{x})^{[t+n]} \det C$. Also
$(\ul{x})^{[t+n]} \det C \subset (\ul{y})^{[l]}$. Therefore $(r_1-r_2)(y_1 \cdots y_d)^{dl}(x_1 \cdots x_d)^t \det C \in (\ul{y})^{[dl+l]}$.
Hence $$(r_1-r_2) \det B (y_1 \cdots y_d)^{dl+l-s} \in (\ul{y})^{[dl+l]}$$ and thus $\ds{[\frac{r_1 \det B}{y_1^{s} \cdots y_d^{s}}]}=\ds{[\frac{r_2 \det B}{y_1^{s} \cdots y_d^{s}}]}$. So $\Phi$ is a well defined map. \qed

\s

We are now ready to examine the relation between the statements in Remark~\ref{equivalences}. We will first assume that $\ul{y}$ form a system of parameters
and determine under which conditions the maps are injective.
The following proposition is proved in \cite[5.1.17]{St}. We include a short proof for the convenience of the reader.

\begin{proposition}{\label{sop->lim inj}}
Let $R$ be a Noetherian local ring, $d=\dim R$ and let $\ul{x}=x_1, \ldots, x_d$ be a system of parameters. Let $\ul{y}=y_1, \ldots, y_d$ be a sequence of elements such that $(\ul{y}) \stackrel{A}{\subset}  (\ul{x})$.  If $\ul{y}$ forms a system of parameters then the map $R/{(\ul{x})^{\rm lim}} \overset{{\cdot \det A}}{\longrightarrow}
 R/{(\ul{y})^{\rm lim}}$ is injective.
\end{proposition}

\proof
Let $r \in R$ such that $r \cdot \det A \in (\ul{y})^{\rm lim}$. Then there exists a positive integer $s$ such that $r \det A (y_1 \cdots y_d)^{s-1} \in (\ul{y})^{[s]}$.
Since $\ul{y}$ is a system of parameters then there exists a matrix $B$ and a positive integer $t$ such that $(\ul{x})^{[t]}=B (\ul{y})^{[s]}$. Hence we have the following inclusions:
\begin{align*}
&(\ul{x})^{[t]} \overset{B}{\subset} (\ul{y})^{[s]} \overset{C}{\subset} (\ul{y}) \overset{A}{\subset} (\ul{x})\\
&(\ul{x})^{[t]} \overset{D}{\subset} (\ul{x}),
\end{align*}
where $C$ is the diagonal matrix with entries $y_i^{s-1}$ and $D$ is the diagonal matrix with entries $x_i^{t-1}$.
Let $E=ACB$. By Corollary~\ref{det-cor} we obtain $(x_1 \cdots x_d)^{td}(\det E -\det D) \in (\ul{x})^{[td+t]}$. Hence $$(x_1 \cdots x_d)^{td}\det B(y_1 \cdots y_d)^{s-1}\det A-(x_1 \cdots x_d)^{td}(x_1 \cdots x_d)^{t-1} \in (\ul{x})^{[td+t]}$$
and thus  $r(x_1 \cdots x_d)^{td}\det B(y_1 \cdots y_d)^{s-1}\det A-r(x_1 \cdots x_d)^{td}(x_1 \cdots x_d)^{t-1} \in (\ul{x})^{[td+t]}$. Again, $$r(x_1 \cdots x_d)^{td}\det B(y_1 \cdots y_d)^{s-1}\det A \in (x_1 \cdots x_d)^{td}\det B (\ul{y})^{[s]} \subset (x_1 \cdots x_d)^{td}(\ul{x})^{[t]} \subset (\ul{x})^{[td+t]}.$$ Therefore $r(x_1 \cdots x_d)^{td}(x_1 \cdots x_d)^{t-1} \in (\ul{x})^{[td+t]}$ and thus
$r \in (\ul{x})^{\rm lim}$. Hence the map is injective. \qed
\s

\begin{remark}\label{Same Phi map}
Let $R$ be a Noetherian local ring, $d=\dim R$ and let $\ul{x}=x_1, \ldots, x_d$ and $\ul{y}=y_1, \ldots, y_d$ be sequences of elements such that $(\ul{y}) \stackrel{A}{\subset}  (\ul{x})$. Then for any   positive integer $n$ there exists a positive integer $s$ such that  $(\ul{y})^{[s]} \subset (\ul{x})^{[n]}$. Then the map $H_{(\ul{x})}^{d}(R) \longrightarrow H_{(\ul{y})}^{d}(R)$ is the same as the map $H_{(\ul{x})^{[n]}}^{d}(R) \longrightarrow H_{(\ul{y})^{[s]}}^{d}(R)$.
\end{remark}

Similarly, we obtain:

\begin{proposition}{\label{sop->lc inj}}
Let $R$ be a Noetherian local ring, $d=\dim R$ and let $\ul{x}=x_1, \ldots, x_d$ be a system of parameters. Let $\ul{y}=y_1, \ldots, y_d$ be a sequence of elements such that $(\ul{y}) \stackrel{A}{\subset}  (\ul{x})$.  If $\ul{y}$ forms a system of parameters then the map 
$H_{(\ul{x})}^{d}(R) \longrightarrow H_{(\ul{y})}^{d}(R)$ defined as in Discussion~\ref{Phi Disc}
is injective.
\end{proposition}

\proof This follows immediately from Proposition~\ref{sop->lim inj}: suppose that  $\Phi(\ds{[\frac{r}{x_1^{n} \cdots x_d^{n}}]} = 0$.
Set $u = \ds{[\frac{1}{x_1^{n} \cdots x_d^{n}}]}$. We recall
how $\Phi$ is defined: 
since $(\ul{y}) \subset (\ul{x})$ there exists a matrix $B$ and a positive integer $s$ such that
$(\ul{y})^{[s]}=B(\ul{x})^{[n]}$. Applying Proposition~\ref{sop->lim inj} to the two systems of
parameters $(\ul{y})^{[s]} \overset{B}{\subset} (\ul{x})^{[n]}$ yields that $ru = 0$. 

To see this notice that $H_{(\ul{x})}^{d}(R)=H_{(\ul{x})^{[n]}}^{d}(R)= \underset{n \in \mathbb{N}}{\bigcup} R/((\ul{x})^{[n]})^{\rm lim}$. Let $s(n)$ be the integer that corresponds to $n$ as above. Applying Proposition~\ref{sop->lim inj} and Remark~\ref{Same Phi map} to the two systems of
parameters $(\ul{y})^{[s]} \overset{B}{\subset} (\ul{x})^{[n]}$ we have:
\[H_{(\ul{x})}^{d}(R)=H_{(\ul{x})^{[n]}}^{d}(R)= \underset{n \in \mathbb{N}}{\bigcup} R/((\ul{x})^{[n]})^{\rm lim}
\hookrightarrow \underset{s(n) \in \mathbb{N}}{\bigcup} R/((\ul{y})^{[s(n)]})^{\rm lim} \hookrightarrow H_{(\ul{y})^{[s]}}^{d}(R)=H_{(\ul{y})}^{d}(R).
\]

\qed

\s

\section{Positive Characteristic}
\bigskip

In this section we are concerned with the maps~(3) and (4) as in Remark~\ref{equivalences}. Recall that * denotes the tight closure.

\begin{proposition}\label{*map}
Let $R$ be a Noetherian local ring of positive characteristic $p$, $d=\dim R$ and let $\ul{x}=x_1, \ldots, x_d$ and $\ul{y}=y_1, \ldots, y_d$ be sequences of elements such that $(\ul{y}) \stackrel{A}{\subset}  (\ul{x})$. Then the map $R/{(\ul{x})^{*}} \overset{{\cdot \det A}}{\longrightarrow} R/{(\ul{y})^{*}}$ is well-defined.
\end{proposition}

\proof
Let $z \in (\ul{x})^{*}$. Then there exists an element $c \in R^{0}$ such that $cz^{q} \in (\ul{x})^{[q]}$ for every $q=p^e$ large.
Since $(\ul{y}) \overset{A}{\subset} (\ul{x})$ then $(\ul{y})^{[q]} \overset{A^{[q]}}{\subset} (\ul{x})^{[q]}$, where $A^{[q]}=(a_{ij}^{q})$.
Notice that $(\det {A})^{q}=\det (A^{[q]})$. Therefore $c(\det {A})^{q}z^{q} \in (\ul{y})^{[q]}$ and thus $(\det A) z \in (\ul{y})^{*}$. \qed

\s

The following Lemma is well-known in the Noetherian case, by linkage theory. We need it in the context of a non-Noetherian ring, however. The proof we give is not the standard one.

\begin{lemma}\label{linkage}
Let $R$ be a commutative ring with identity and let $I$ be an $R$-ideal. Suppose that $\ul{y}=y_1, \ldots, y_d$ is a regular sequence and $(\ul{y}) \subset I$. Let $0 \lr R^{n} \lr \cdots \lr R \lr R/I \lr 0$ be a free resolution of $R/I$. We further assume that the following diagram has exact rows:
\begin{diagram}[width=2em]
&0 &\rTo &R^{n} &\rTo &\cdots &\rTo &R &\rTo &R/I &\rTo &0\\
&& &\uTo^{\left(\begin{array}{c}f_1 \\ \vdots\\ f_n\\ \end{array} \right )} &&&& \parallel & & \uTo &&\\
&0 &\rTo &R &\rTo^{\ds{(y_1 \ldots y_d)}} &\cdots &\rTo &R &\rTo &R/(\ul{y}) &\rTo &0\\
\end{diagram}

Then $(\ul{y}):I=(\ul{y})+(\ul{f})$, where $\ul{f}=f_1, \ldots, f_n$.
\end{lemma}

\proof
First notice that $R$ need not be Noetherian. To see that $(\ul{y})+(\ul{f}) \subset (\ul{y}):I$, let $i \in I$ and extend the natural  map  $R/I \stackrel{\cdot i} {\lr} R/I$ given by multiplication by $i$ to a map from the free resolution of $R/I$ to the free resolution of $R/I$ as in the diagram below. The composition of the map $R/(\ul{y}) \lr R/I$ with the map $R/I \stackrel{\cdot i} {\lr} R/I$ is the zero map. We begin defining a homotopy by taking  $\sd_0=0$:

\begin{diagram}[width=2em]
&0 &\rTo &R^{n} &\rTo &\cdots &\rTo &R &\rTo &R/I &\rTo &0\\
&& &\uTo^{\cdot i} &&&& \uTo^{\cdot i} &\luTo(2,4)\;\;\;\;\;{\sd_{0}=0} & \uTo^{\cdot i}  &&\\
&0 &\rTo &R^{n} &\rTo &\cdots &\rTo &R &\rTo &R/I &\rTo &0\\
&& & \uTo^{\left(\begin{array}{c}f_1 \\ \vdots\\ f_n\\ \end{array} \right )} &&&&\parallel & & \uTo&&\\
&0 &\rTo &R &\rTo^{\ds{(y_1 \ldots y_d )}} &\cdots &\rTo &R &\rTo &R/(\ul{y}) &\rTo &0\\
\end{diagram}

By Lemma~\ref{homotopyfrombelow} we may extend the homotopy to obtain $i(f_1, \ldots, f_d) \subset (y_1, \ldots, y_d)$.

For the other inclusion let $s \in (\ul{y}):I$ and consider the map $R/I \stackrel{\cdot s}{\lr} R/(\ul{y})$ given by multiplication by $s$. Extending this map  we obtain the following diagram:

\begin{diagram}[width=2em]
&0 &\rTo &R &\rTo &\cdots &\rTo &R &\rTo &R/(\ul{y}) &\rTo &0\\
&& &\uTo^{\left(\begin{array}{c}g_1  \ldots g_n\\ \end{array} \right )} &&&& \uTo^{\cdot s} & & \uTo^{\cdot s} &&\\
&0 &\rTo &R^{n} &\rTo &\cdots &\rTo &R &\rTo &R/I &\rTo &0\\
&& &\uTo^{\left(\begin{array}{c}f_1 \\ \vdots\\ f_n\\ \end{array} \right )} &&&& \parallel & & \uTo &&\\
&0 &\rTo &R &\rTo &\cdots &\rTo &R &\rTo &R/(\ul{y}) &\rTo &0\\
\end{diagram}

for some elements $g_1, \ldots, g_n \in R$.

On the other hand we consider the map $R/(\ul{y}) \stackrel{\cdot s}{\lr} R/(\ul{y})$ given by multiplication by $s$. Extending this map we obtain the following diagram:

\begin{diagram}[width=2em]
&0 &\rTo &R &\rTo &\cdots &\rTo &R &\rTo &R/(\ul{y}) &\rTo &0\\
&& &\uTo^{\cdot s} &&&& \uTo^{\cdot s} & & \uTo^{\cdot s} &&\\
&0 &\rTo &R &\rTo &\cdots &\rTo &R &\rTo &R/(\ul{y}) &\rTo &0\\
\end{diagram}

Combining the two diagrams we can define a homotopy by taking $\sd_0=0$ and by Lemma~\ref{homotopyfrombelow} we can extend the homotopy:

\begin{diagram}[width=2em]
&0 &\rTo &R &\rTo &\cdots &\rTo &R &\rTo &R/(\ul{y}) &\rTo &0\\
&& &\uTo^{\left(\begin{array}{c}f_1 \\ \vdots\\ f_n\\ \end{array} \right )\left(\begin{array}{c}g_1  \ldots g_n\\ \end{array} \right )-\cdot s} & \luTo(2,2)\;{\sd_{m}} &&& \uTo^{0}& \luTo(2,2)\; \; \; \; \; \; \; \; \; \; \; \; \; \; \; \; \;\;{\sd_{0}=0} & \uTo &&\\
&0 &\rTo &R &\rTo^{(y_1, \ldots, y_d) \; \; \; \; \; \; \; \;} &\cdots &\rTo &R &\rTo &R/(\ul{y}) &\rTo &0\\
\end{diagram}

Therefore $s-\sum_{i=1}^{n} f_ig_i \in (y_1, \ldots, y_d)$ and  $s \in (\ul{y})+(\ul{f})$.

\qed

\s

\begin{proposition}{\label{sop->*clos}}
Let $R$ be an excellent local ring of positive characteristic $p$, $d=\dim R$, and 
let $\ul{x}=x_1, \ldots, x_d$ be a system of parameters. Let $\ul{y}=y_1, \ldots, y_d$ be a sequence of elements such that $(\ul{y}) \stackrel{A}{\subset}  (\ul{x})$. 
If $\ul{y}$ forms a system of parameters then the map $R/{(\ul{x})^{*}} \overset{{\cdot \det A}}{\longrightarrow} R/{(\ul{y})^{*}}$ is injective.
\end{proposition}

\proof
We can first pass to the completion and thus assume that $R$ is a complete. By Proposition~\ref{*map} the map
$R/{(\ul{x})^{*}} \overset{{\cdot \det A}}{\longrightarrow} R/{(\ul{y})^{*}}$ is well defined and thus
$(\ul{x})^{*} \subset (\ul{y})^{*}:\det A$.  Let $r \in (\ul{y})^{*}:\det A$. To show that $r \in (\ul{x})^{*}$ it is enough
to show that the image of $r$ in $R/P$ is in the tight closure of $((\ul{x})+P)/P$ for every minimal prime $P$
of $R$ by \cite[Proposition~4.1]{HoHu1}. Thus we may assume that $R$ is a domain.

Since $r\det A \in (\ul{y})^{*}$ and $(\ul{y})^{*} \subset (\ul{y})R^{+}$ then $r\det A \in (\ul{y})R^{+}$. 
In $R^{+}$ both $\ul{x}$ and $\ul{y}$ form a regular sequence. Therefore $(\ul{y})R^{+}:\det A=(\ul{x})R^{+}$ since
$R^{+}$ is Cohen-Macaulay or by Lemma~\ref{linkage}. Hence $r \in (\ul{x})R^{+} \cap R=(\ul{x})^{*}$, where the equality follows by \cite[Theorem~5.1]{Sm}. \qed

\s

The fact that parameters form a regular sequence in $R^+$ immediately gives injectivity on local cohomology:

\begin{proposition}{\label{sop->lc+ inj}}
Let $R$ be a Noetherian local domain of positive characteristic $p$, $d=\dim R$, and 
let $\ul{x}=x_1, \ldots, x_d$ be a system of parameters. Let $\ul{y}=y_1, \ldots, y_d$ be a sequence of elements such that $(\ul{y}) \stackrel{A}{\subset}  (\ul{x})$. 
If $\ul{y}$ forms a system of parameters then the map $H_{(\ul{x})}^{d}(R^+) \longrightarrow H_{(\ul{y})}^{d}(R^+)$
is injective.
\end{proposition}

\proof This follows as in  Proposition~\ref{sop->lc inj}: suppose that  $\Phi(\ds{[\frac{r}{x_1^{n} \cdots x_d^{n}}]} = 0$.
Set $u = \ds{[\frac{1}{x_1^{n} \cdots x_d^{n}}]}$. We recall
how $\Phi$ is defined:
since $(\ul{y}) \subset (\ul{x})$ there exists a matrix $B$ and a positive integer $s$ such that
$(\ul{y})^{[s]}=B(\ul{x})^{[n]}$.
Since parameters form a regular sequence in $R^+$,  $(\ul{y})^{[s]}:\det B =  (\ul{x})^{[n]}$, and then
it follows that  $ru = 0$. \qed

\bigskip

\section{The one dimensional case}

In this section we prove two main results, that the injectivity of either map
$R/(x) \overset{{\cdot u}}{\longrightarrow} R/(y)$ or $ R/{(x)^{\rm lim}} \stackrel{ \cdot u}{\lr} R/{(y)^{\rm lim}}$
forces $y$ to be a parameter.  It is interesting that even this simple case is not obvious.

\s

\begin{theorem}\label{1-dim DR general}
Let $R$ be a 1-dimensional Noetherian local ring. Let $x$ be a parameter and let $y=ux$.
 If the map $R/(x) \overset{{\cdot u}}{\longrightarrow} R/(y)$ is injective then $y$ is a parameter.
\end{theorem}

\proof
Suppose that $x$ is a parameter and that the map $R/(x) \overset{{\cdot u}}{\longrightarrow} R/(y)$ is injective. Since $x$ is a parameter, 
$0:x$ has finite length. Consider the following exact sequence:
\[
0 \lr (0:(x,u)) \lr 0:x \stackrel{ \cdot u}{\lr} 0:x \lr \ds{\frac{0:x}{u(0:x)}} \lr 0 \ .
\]

Computing length we obtain that  $\la (\ds{\frac{0:x}{u(0:x)}})=\la ( (0:(x,u)))$. Suppose that $u$ is not a parameter. Then
there exists a minimal prime $P$ such that $u \in P$. Hence $0:(0:u) \subset P$ and in particular $ \dim ((0:u)) >0$, 
and consequently $e(x;(0:u)) >0$.  On the other hand, $e(x;(0:u))=\la(\ds{\frac{0:u}{x(0:u)}}) - \la(0:_{(0:u)} x) >0$. 
 Since $0:_{(0:u)} x=0:(x,u)$ then by the above computations we obtain $\la(\ds{\frac{0:u}{x(0:u)}}) > \la (\ds{\frac{0:x}{u(0:x)}})$.

Now since the map is injective we have $0:u \subset (x)$ and thus $0:u=x(0:xu)$. Hence
\[\ds{\frac{0:u}{x(0:u)}=\frac{x(0:xu)}{x(0:u)} \simeq \frac{0:xu}{(0:u)+(0:x)} \simeq \frac{u(0:xu)}{u(0:x)} \subset \frac{0:x}{u(0:x)}}\]

Therefore $\la(\ds{\frac{0:x}{u(0:x)}}) \geq \la(\ds{\frac{u(0:xu)}{u(0:x)}})=\la(\ds{\frac{0:u}{x(0:u)}})$, which is a contradiction. 
Thus $u$ must be a parameter and in conclusion $y$ is also a parameter. \qed

\begin{corollary}\label{ulrich-question}
Let $R$ be a 1-dimensional Noetherian local ring. Let $x$ be a parameter and let $y=ux$. 
The map $R/(x) \overset{{\cdot u}}{\longrightarrow} R/(y)$ is injective if and only if the map $R/(u) \overset{{\cdot x}}{\longrightarrow} R/(y)$ is injective.
\end{corollary}

\proof
Notice that it suffices to show one direction. Suppose that $R/(x) \overset{{\cdot u}}{\longrightarrow} R/(y)$ is injective. Then by 
Theorem~\ref{1-dim DR general} $u$ is also a parameter.
By the proof of Theorem~\ref{1-dim DR general} we see that by symmetry one has
\[
\ds{\frac{0:u}{x(0:u)}\simeq \frac{x(0:xu)}{x(0:u)} \simeq \frac{0:xu}{(0:x)+(0:u)} \simeq \frac{u(0:xu)}{u(0:x)} \subset \frac{0:x}{u(0:x)}}\]
But since $\la(\ds{\frac{0:u}{x(0:u)}})=\la(\ds{\frac{0:x}{u(0:x)}})$, then $\ds{\frac{u(0:xu)}{u(0:x)} \simeq \frac{0:x}{u(0:x)}}$. 
Therefore $0:x \subset (u)$, i.e. $R/(u) \overset{{\cdot x}}{\longrightarrow} R/(y)$ is injective. \qed

\s

The converse of Theorem~\ref{1-dim DR general} does not hold in general.
The following example was shown to us by R. Heitmann:

\begin{example}\label{ex Heitmann}
Let $R=k[x,u]/(((x+u)u)^3,x(x+u)^2u^2)$, where ${\rm char} \; k=2$. Notice that $R$ is a one-dimensional ring
where $x$ is a parameter. Let $y=x^2$ be also a parameter. Then one can see that $(y):(x)=(x,u^4)$, which means that
the map $R/(x) \stackrel{\cdot x}{\lr} R/(y)$ is not injective.
\end{example}

In fact, in Section 6 we will prove in arbitrary dimension that if the maps are injective for all parameters, then
the ring must be Cohen-Macaulay.

\s

\begin{theorem}\label{1-dim iff lim clo}
Let $R$ be a 1-dimensional Noetherian local ring. Let $x$ be a parameter, and
let $y=ux$. Then $y$ is a parameter if and only if the map $ R/{(x)^{\rm lim}} \stackrel{ \cdot u}{\lr} R/{(y)^{\rm lim}}$ is injective.

\end{theorem}

\proof 
First notice that the forward direction follows from Proposition~\ref{sop->lim inj}. So we assume that the map
$ R/{(x)^{\rm lim}} \stackrel{ \cdot u}{\lr} R/{(y)^{\rm lim}}$ is injective.
Note that $(x)^{\rm lim} = (x^{n+1}):x^n$ for all large $n$, which in turn is equal to $(x) + (0:x^n)$ for large $n$.
Similarly, $(y)^{\rm lim} = (y) + (0:y^n)$ for large enough $n$. We can analyze these ideals using a primary decomposition
of $(0)$. Let $0 = q_1\cap q_2\cap \ldots\cap q_s\cap J$, where $\sqrt{q_i} = p_i$ are distinct minimal primes, and
$\sqrt{J} = \m$. We wish to prove that $y$ is a parameter. Assume not. Then $y$ is in at least one minimal prime. Let
$y\in p_1\cap\ldots\cap p_t$, and $y\notin p_{t+1}\cup \ldots\cup p_s$. We see that  $(x) + (0:x^n) = (x) + q_1\cap q_2\cap \ldots\cap q_s$,
and $(y) + (0:y^n) = (y) +  q_{t+1}\cap \ldots\cap q_s$. The assumption on the injectivity of the map then becomes that
$((y) +  q_{t+1}\cap \ldots\cap q_s):u\subset (x) +  q_1\cap q_2\cap \ldots\cap q_s$. Since $u$ is not in $ p_{t+1}\cup \ldots\cup p_s$,
it follows that  $((y) +  q_{t+1}\cap \ldots\cap q_s):u = (x) +  q_{t+1}\cap \ldots\cap q_s\subset (x) + q_1\cap q_2\cap \ldots\cap q_s$.
In particular, $ q_{t+1}\cap \ldots\cap q_s\subset x( q_{t+1}\cap \ldots\cap q_s) + q_1\cap q_2\cap \ldots\cap q_s$, which by Nakayama's
lemma shows that $t = 0$. It follows that $y$ is a parameter.  \qed

\s

\begin{discussion}{\rm It is natural to believe that the one-dimensional case above would at least prove the higher dimensional
case when the matrix $A$ is a diagonal matrix, i.e., the case in which $x_1,\ldots,x_d$ are parameters, and
$y_i = u_ix_i$. In this case our assumption would be that the map
$R/(\ul{x})\overset{u_1\cdots u_d}\to R/(\ul{y})$ is injective. One is tempted to break this up into several
maps by changing one $x_i$ at a time. Although this works for one step, it seems to break down even in dimension $2$, and
gives a good idea of the difficulty in extending to higher dimensions. 
}\end{discussion}
\s
\section{Higher Dimensions}
We would like to extend Theorem~\ref{1-dim iff lim clo} to higher dimensions. We are able to do so in the case the ring
is analytically irreducible, as well as the case in which  the system of parameters $\ul{x}$ lies deep inside the maximal ideal.

\begin{proposition}\label{Hart-Lich iff}
Let $R$ be an analytically irreducible local ring of equicharacteristic and  $d=\dim R$, or
of mixed characteristic and of dimension at most $3$. let $\ul{x}=x_1, \ldots, x_d$ be a system of parameters. Let $\ul{y}=y_1, \ldots, y_d$ be a sequence of elements such that $(\ul{y}) \stackrel{A}{\subset}  (\ul{x})$.  Then $\ul{y}$ forms a system of parameters if and only if the map
$R/{(\ul{x})^{\rm lim}} \overset{{\cdot \det A}}{\longrightarrow} R/{(\ul{y})^{\rm lim}}$ is injective.
\end{proposition}

\proof We may assume that $R$ is a complete domain. The forward direction is already covered in
Proposition~\ref{sop->lim inj}.

Suppose that the map $R/{(\ul{x})^{\rm lim}} \overset{{\cdot \det A}}{\longrightarrow} R/{(\ul{y})^{\rm lim}}$ is injective.
Since $\ul{x}$ is a system of parameters then $(\ul{x})^{\rm lim} \neq R$; this is the monomial conjecture of Hochster, which
is true in either equicharacteristic (\cite{Ho}) or in dimension at most three \cite{He}. As the map
$R/{(\ul{x})^{\rm lim}} \overset{{\cdot \det A}}{\longrightarrow} R/{(\ul{y})^{\rm lim}}$ is injective then
$R/{(\ul{y})^{\rm lim}} \neq 0$. Also, notice that $R/{(\ul{y})^{\rm lim}}$ injects naturally into
$H_{(\ul{y})}^{d}(R)$ and thus $H_{(\ul{y})}^{d}(R) \neq 0$. Therefore by Hartshorne-Lictenbaum Vanishing Theorem
there exists a minimal prime $P$ of $R$ with $\dim \; R/P =\dim \; R$ such that $\dim \; (R/((\ul{y}) +P))=0$.
Since $R$ is a domain it follows that $\dim \; R/(\ul{y})=0$, and thus $\ul{y}$ is a system of parameters. \qed

\s
If $R$ is not analytically irreducible, then unfortunately we cannot reach the same conclusion. See Example~\ref{highpower}
in the last section. However, we can say something:

\begin{theorem}\label{lim inj iff Hd inj}
Let $(R,\m)$ be a Noetherian local ring, $d=\dim R$, and let $\ul{x}=x_1, \ldots, x_d$ be a system of parameters. Let $\ul{y}=y_1, \ldots, y_d$ be a sequence of elements such that $(\ul{y}) \stackrel{A}{\subset}  (\ul{x})$. 
There exists an integer $\ell$ such that if $(\ul{x}) \subset \m^{\ell}$ then the map
$R/{(\ul{x})^{\rm lim}} \overset{{\cdot \det A}}{\longrightarrow} R/{(\ul{y})^{\rm lim}}$ is injective if and only if
the map $H_{(\ul{x})}^{d}(R) \lr H_{(\ul{y})}^{d}(R)$ as in Discussion~\ref{Phi Disc} is injective.

\end{theorem}

\proof

By \cite[Lemma~3.12]{GoSa} there exists an integer $\ell$ such that if $\ul{x}$ is system of
parameters in $\hat{R}$ such that $(\ul{x}) \subset \m^{\ell}\hat{R}$ then the map ${\rm Soc}(\hat{R}/(\ul{x})) \lr
{\rm Soc}( H_{(\ul{x})}^{d}(\hat{R}))$ is surjective. Notice that $\ul{x}$ is a system of parameters in $\hat{R}$ if
and only if $\ul{x}$ is a system of parameters in $R$. Also the map $H_{(\ul{x})}^{d}(R) \lr H_{(\ul{y})}^{d}(R)$ is injective if
and only if $H_{(\ul{x})}^{d}(\hat{R}) \lr H_{(\ul{y})}^{d}(\hat{R})$ is injective. Therefore without loss of generality we may assume that $R$ is complete.

Suppose that $(\ul{x}) \subset \m^{\ell}$. Notice that $R/{(\ul{x})^{\rm lim}}$ injects naturally into $H_{(\ul{x})}^{d}(R)$ and similarly $R/{(\ul{y})^{\rm lim}}$ injects naturally into $H_{(\ul{y})}^{d}(R)$.
By Proposition~\ref{Hmap} there exists a homomorphism $ \Phi: H_{(\ul{x})}^{d}(R) \lr H_{(\ul{y})}^{d}(R)$ and moreover the following diagrams commute:

\begin{diagram}[width=2em]
0 &\rTo& R/{(\ul{x})^{\rm lim}} & \rTo^{\alpha} & H_{(\ul{x})}^{d}(R) \\
&& \dTo^{\beta=\cdot \det A} & & \dTo^{\Phi} \\
0 &\rTo& R/{(\ul{y})^{\rm lim}} & \rTo^{\gamma} & H_{(\ul{y})}^{d}(R) \\
\end{diagram}

and
\begin{diagram}[width=2em]
0 &\rTo &\ds{\frac{(\ul{x})^{\rm lim}}{(\ul{x})}} &\rTo &\ds{\frac{R}{(\ul{x})}} &\rTo^{\phi_{\ul{x}}} &H_{(\ul{x})}^{d}(R) \\
&&& &\dTo^{\gamma= \cdot \det A} & & \dTo^{\Phi} \\
0 &\rTo &\ds{\frac{(\ul{y})^{\rm lim}}{(\ul{y})}} &\rTo &\ds{\frac{R}{(\ul{y})}} &\rTo^{\phi_{\ul{y}}} &H_{(\ul{y})}^{d}(R) \\
\end{diagram}

Notice that if $\Phi$ is injective then clearly $\beta$ is injective. Suppose that the map $R/{(\ul{x})^{\rm lim}} \overset{{\cdot \det A}}{\longrightarrow} R/{(\ul{y})^{\rm lim}}$ is injective. Suppose that $\Phi$ is not injective. Since $ H_{(\ul{x})}^{d}(R)$ is an Artinian module then $ H_{(\ul{x})}^{d}(R)$ is essential over ${\rm Soc}(H_{(\ul{x})}^{d}(R))$. Hence there exists a nonzero element $f \in \ker (\Phi) \cap {\rm Soc}( H_{(\ul{x})}^{d}(R))$. By \cite[Lemma~3.12]{GoSa} the map ${\rm Soc}(R/{(\ul{x})}) \stackrel{\phi_{\ul{x}}}\lr {\rm Soc}( H_{(\ul{x})}^{d}(R))$ is surjective. Hence there exists $g \in {\rm Soc}(R/{(\ul{x})})$ such that $\phi_{\ul{x}}(g)=f$. But as the second diagram commutes then $\phi_{\ul{y}}(\gamma(g))=\Phi(\phi_{\ul{x}}(f))$, or in other words $\phi_{\ul{y}}((\det A) g)=0$. Since $\ker \phi_{\ul{y}}=\ds{\frac{(\ul{y})^{\rm lim}}{(\ul{y})}}$ then $(\det A) g \in (\ul{y})^{\rm lim}$. Thus since $\beta$ is injective then $g \in (\ul{x})^{\rm lim}$ or in other words $g \in \ker \phi_{\ul{x}}$ and hence $f=0$, which is a contradiction. Therefore, $\Phi$ is also injective.\qed

\s

\begin{theorem}
Let $(R,\m)$ be a Noetherian local ring  and let $d=\dim R$. There exists an integer $\ell$ with the following property:
whenever $\ul{x}=x_{1}, \ldots, x_d$ is a system of parameters with $(\ul{x}) \subset \m^{\ell}$ and
$\ul{y}=y_1, \ldots, y_d$ a sequence of elements such that $(\ul{y}) \stackrel{A}{\subset} (\ul{x})$, then
$\ul{y}$ forms a system of parameters if and only if the map $H_{(\ul{x})}^{d}(R) \lr H_{(\ul{y})}^{d}(R)$ as in  Discussion~\ref{Phi Disc} is injective.

\end{theorem}

\proof
As in the proof of Theorem~\ref{lim inj iff Hd inj} we may pass to the completion $\hat{R}$ and assume that $R$ is complete. Let $\ell$ be the integer as in \cite[Lemma~3.12]{GoSa} such that $(\ul{x}) \subset \m^{\ell}$ and the map ${\rm Soc}(R/(\ul{x})) \lr {\rm Soc}( H_{(\ul{x})}^{d}(R))$ is surjective.

First we assume that $\ul{y}$ is a system of parameters. Then by Proposition~\ref{sop->lim inj} the map
$R/{(\ul{x})^{\rm lim}} \overset{{\cdot \det A}}{\longrightarrow} R/{(\ul{y})^{\rm lim}}$ is injective
and thus by Theorem~\ref{lim inj iff Hd inj} the map $H_{(\ul{x})}^{d}(R) \lr H_{(\ul{y})}^{d}(R)$ as in Proposition~\ref{Hmap} is injective.

Suppose the map $H_{(\ul{x})}^{d}(R) \lr H_{(\ul{y})}^{d}(R)$ as in Proposition~\ref{Hmap} is injective. To show that $\ul{y}$
is a system of parameters,  it suffices to show that $\ul{y}$ is a system of parameters in $\ol{R}=R/P$ for every minimal prime $P$ of maximal dimension. Suppose that it is not true. Then there exists a minimal prime prime $P$ of maximal dimension such that $H_{(\ul{y})}^{d}(\ol{R}) =0$, where $\ol{R}=R/P$. Let $P=0:c$ for some $c\in R$. We have the following short exact sequence:
\[
0 \lr R/(0:c) \lr R \stackrel{\pi}{\lr} R/(c) \lr 0\;.
\]
Since $H_{(\ul{y})}^{d}(\ol{R}) =0$ we obtain the sequence $H_{(\ul{x})}^{d}(R) \stackrel{\Phi}\hookrightarrow H_{(\ul{y})}^{d}(R) \stackrel{\pi}{\hookrightarrow} H_{(\ul{y})}^{d}(R/(c)) \lr 0$.

We claim that $c H_{(\ul{x})}^{d}(R)=0$. Let $f \in H_{(\ul{x})}^{d}(R)$. Then $\pi(\Phi(c f))=\pi(c \Phi(f))=c\pi(\Phi(f))$.
But $c H_{(\ul{y})}^{d}(R/c)=0$ so $c\pi(\Phi(f))=0$. On the other hand though both $\pi$ and $\Phi$ are injective and hence
$cf=0$, which shows the claim.

We note now that since $c H_{(\ul{x})}^{d}(R)=0$ it follows that  $c (H_{(\ul{x})}^{d}(R))^{\vee} = 0$. But
the dual of the top local cohomology is $\omega_{R}$, the canonical module of $R$. Thus $c \in {\rm ann} \omega_R$, which is a contradiction,
since by \cite[Remark~2.2]{HoHu2} ${\rm ann} \omega_{R}=\{r \in R \; | \; \dim \; R/(0:r) <d \}$.

\qed

\s

\begin{corollary} \label{lim inj iff par in m^l}
Let $(R,\m)$ be a Noetherian local ring  and let $d=\dim R$. There exists an integer $\ell$ with the following property:
whenever $\ul{x}=x_{1}, \ldots, x_d$ is a system of parameters with $(\ul{x}) \subset \m^{\ell}$ and
$\ul{y}=y_1, \ldots, y_d$ a sequence of elements such that $(\ul{y}) \stackrel{A}{\subset} (\ul{x})$, then 
 $\ul{y}$ forms a system of parameters if and only if the map $R/{(\ul{x})^{\rm lim}} \overset{{\cdot \det A}}{\longrightarrow} R/{(\ul{y})^{\rm lim}}$ is injective.

\end{corollary}

\proof Simply combine the last two theorems. \qed

\s

\section{Examples, Extensions, and Questions}

\s
In this section we conclude with several examples which show that the hypotheses of several of the theorems
are necessary, consider some extensions of our results and list some open questions.

The first example shows that the condition that the parameters are in a deep enough power of the maximal ideal
is necessary in Corollary~\ref{lim inj iff par in m^l}.

\begin{example}\label{highpower}
Let $R=\ds{\frac{k[a,b,c,d]}{(a,b) \cap (c,d)}}$, where $k$ is a field. Notice that $R$ is an equidimensional ring. 
Let $x_1=a+c, x_2=b+d$ and $y_1=a^2,y_2=b^2$. We claim that $\ul{x}=x_1,x_2$ is a system of parameters, $(\ul{y})=(y_1,y_2)
\subset (\ul{x})$, the map $R/{(\ul{x})^{\rm lim}} \overset{{\cdot \det A}}{\longrightarrow} R/{(\ul{y})^{\rm lim}}$
is injective and $\ul{y}$ is not a system of parameters.

First note that it is straight forward to see that $\ul{x}=x_1,x_2$ is a system of parameters, $(\ul{y})=(y_1,y_2) \subset
(\ul{x})$ and $\ul{y}$ is not a system of parameters. We now claim that $(\ul{x})^{\rm lim}=\m$. Notice that $\m (x_1x_2) \subset
(x_1^2,x_2^2)$ and thus $\m \subset (\ul{x})^{\rm lim}$. Since $(\ul{x})^{\rm lim} \neq R$ then $(\ul{x})^{\rm lim}=\m$.

We also claim that $(\ul{y})^{\rm lim}=(a^2,b^2,c,d)$. Notice that $(a^2,b^2,c,d)(a^2b^2) \subset (a^4,b^4)$. On the
other hand suppose $z (a^{2n}b^{2n}) \subset (a^{2n+2},b^{2n+2})$. We consider the ring $\overline{R}=R/(c,d) \simeq k[a,b]$.
Let $^{^{\tratto}}$ denote the images in $\overline{R}$. Then $\ol{z} \subset (a^2,b^2) \ol{R}$. Thus
$z \in (a^2,b^2,c,d)$ and therefore $(a^2,b^2,c,d)=(\ul{y})^{\rm lim}$.

We now consider the map $R/{(\ul{x})^{\rm lim}} \overset{{\cdot \det A}}{\longrightarrow} R/{(\ul{y})^{\rm lim}}$, where
$\det A=ab$. Notice that $1$ maps to $ab$ and the map is clearly injective.
\end{example}

The next example deals with the question of whether or not the injectivity of the map in the original situation
of Dutta and Roberts actually forces the ring to be Cohen-Macaulay. The answer is no, as the first example shows,
but is yes if one requires injectivity for all systems of parameters.

\begin{example}
Let $R=k[x,z]/(x^2z,z^2)$ be a one-dimensional ring. Then $R$ is not Cohen-Macaulay.  However, $x$ and $y=x^2$ are parameters, 
and $x^2:_Rx=(x)$. 

\end{example}
\s

\begin{proposition}\label{all sop->CM}
Let $R$ be a Noetherian local ring  and let $d=\dim R$. Suppose that for every system
of parameters $\ul{x}=x_{1}, \ldots, x_d$ and every positive integer $t$
the map $R/(\ul{x}) \stackrel{(x_1 \cdots x_d)^{t-1}}{\lr} R/(\ul{x})^{[t]}$ is injective. Then $R$ is Cohen-Macaulay.
\end{proposition}

\proof
Notice that $(\ul{x})^{\rm lim}=(\ul{x})^{[t]}:(x_1 \cdots x_d)^{t-1}$ for some $t \in  \mathbb{N}$. 
Since $\ul{x}$ is a system of parameters then so is $\ul{x}^{[t]}$. Therefore by assumption the map
$R/(\ul{x}) \stackrel{\det D}{\lr} R/(\ul{x})^{[t]}$ is injective, where now $D$ is the diagonal matrix
with diagonal entries $x_i^{t-1}$. Then $(\ul{x})=(\ul{x})^{[t]}:\det D=(\ul{x})^{[t]}:(x_1 \cdots x_d)^{t-1}$ since
the map is injective. Hence $(\ul{x})=(\ul{x})^{\rm lim}$ and by \cite[Proposition~2.3]{MRS}(see also \cite[Theorem~5.2.3]{St}) $R$ is then Cohen-Macaulay. \qed

\s

We close with two additional questions which are suggested by the work in this paper.

\s
\begin{question}
Let $R$ be a Noetherian local ring. Suppose that $u$ is in some minimal prime $P$ such that
the dimension of $R/P$ is the same as that of $R$. Can $0:u$ ever be in an ideal generated by a system of
parameters?

\end{question}

\s

We proved that the answer is `no' in dimension one, and our calculations strongly suggest the answer is
always `no'.  A related observation is due to Strooker and Simon \cite{SiS}: 
They prove the following: let $R$ be a Gorenstein local ring and let 
$A = R/I$ where $I$ in a nonzero ideal in $R$ consisting of
zero divisors. Set $J = 0:I$.
Then $A$ satisfies
the monomial conjecture if and only if the ideal $J$ is not contained in any parameter ideal of $R$.

Another question was suggested to us by Bernd Ulrich.
\begin{question}
 Our main question can be rephrased to say that
$(\ul{y}):\det A = (\ul{x})$ should imply that $\ul{y}$ are parameters. This statement, in the Cohen-Macaulay case, is equivalent to the dual statement that
$(\ul{y}):(\ul{x}) = (\ul{y},\det A)$. What is the relationship between the two in  general? 
\end{question}

Note that 
in the one-dimensional case, Corollary~\ref{ulrich-question} proves the two are equivalent.
\s

\end{document}